\documentclass{amsart}
\usepackage{amsfonts}
\usepackage{latexsym}
\usepackage{amsmath}
\usepackage{amssymb}
\usepackage{amsthm}
\usepackage{fullpage}
\usepackage{enumerate}
\usepackage{color}
\newtheorem{theorem}{Theorem}[section]
\newtheorem{lemma}[theorem]{Lemma}

\newtheorem{proposition}[theorem]{Proposition}

\newtheorem{remark}[theorem]{Remark}

\newtheorem{question}[theorem]{Question}
\theoremstyle{definition}

\newtheorem{thm}{Theorem}[section] 
\newtheorem{thmy}{Theorem}
\newenvironment{oldtheorem}{\stepcounter{thm}\begin{thmy}}{\end{thmy}}
\newtheorem*{note*}{Note}

\makeatletter

\def\R{{\mathbb R}}
\def\endproof{\begin{flushright}
$ \Box $ \\
\end{flushright}}

\makeatother

\begin{document}


\title{Estimating volume and surface area of a convex body via its
projections or sections}

\medskip

\dedicatory{Dedicated to the memory of our friend  and colleague Joe Diestel, 1943-2017.}

\author {Alexander  Koldobsky,  Christos Saroglou, and  Artem Zvavitch}\thanks{The third named author is supported in part by  the U.S. National Science Foundation Grant DMS-1101636.}

\address{Department of Mathematics, University of Missouri, Columbia, MO 65211, USA} \email{koldobskiya@missouri.edu}

\address{Department of Mathematics, Kent State University,
Kent, OH 44242, USA} \email{csaroglo@kent.edu  \&  christos.saroglou@gmail.com}

\address{Department of Mathematics, Kent State University,
Kent, OH 44242, USA} \email{zvavitch@math.kent.edu}



\begin{abstract}
\footnotesize 


The main goal of this paper is to present a series of  inequalities connecting the  surface area measure of a convex body and surface area measure of its projections and sections. 
We present a solution of a question from \cite{CGG}  regarding the asymptotic behavior of the best constant in a recently proposed reverse Loomis-Whitney inequality. Next we   give a new sufficient condition for the slicing problem to have an affirmative answer, in terms of the least ``outer volume ratio distance'' from the class of intersection bodies of projections of at least proportional dimension of convex bodies. Finally,  we show that certain geometric quantities such as the volume ratio and minimal surface area (after a suitable normalization) are not necessarily close to each other. 
\end{abstract}

\maketitle

\section{Introduction}
\hspace*{1.5em}In the past decades a lot of effort has been put in the study of problems of estimating volumetric quantities of a convex body (i.e. a convex compact set with non-empty interior) in terms of the corresponding functionals of its sections or projections. We refer to the books \cite{Ga, Ko1, KoY, RZ, S} for more information, examples and the history of those problems. Our aim is to continue this investigation by considering a number of  problems of this type. 

In Section 3, we study inequalities involving the size of  projections of $n$-dimensional convex body and study the following problem proposed in  \cite{CGG}: 
\begin{question}\label{qu:lw}
Find the smallest constant $\Lambda_n$, such that the following inequality holds for all convex bodies in $\mathbb{R}^n$:
\begin{equation}\label{eq:lum}
\min_{\{e_1,\dots,e_n\}\in{\mathcal F}_n}\prod_{i=1}^n\big|K|e_i^{\perp}\big|\leq \Lambda_n|K|^{n-1},
\end{equation}
where ${\mathcal F}_n$ denotes the set of all orthonormal basis' in $\mathbb{R}^n$. 
\end{question}
 Inequality  (\ref{eq:lum}) can be viewed as a reverse to the classical Loomis-Whitney \cite{LW} inequality for compact subset $A \subset \mathbb{R}^n$:
 \begin{equation}\label{eq:lw}
\prod_{i=1}^n\big|A|e_i^{\perp}\big|\ge  |A|^{n-1}.
\end{equation}
 The authors in \cite{CGG} asked for the correct asymptotical behavior of $\Lambda_n^{1/n}$.  In section 3  (Theorem \ref{thm-reverse-lw})  we show that $\Lambda_n^{1/n}\leq c\sqrt{n}$, for some absolute constant $c>0$. Moreover, we prove that this estimate is the best possible up to an absolute  constant $c$. We note that other  variants of reverse Loomis-Whitney type inequalities where considered in  \cite{CGaG}.

In Section 4, we  study a number of inequalities which arise in the study of comparison problems. For example, what information on convex bodies we can get  from comparison inequalities for its curvature functions.  We also study the relationship of volume ratio, curvature measure and surface area of convex bodies. For instance we prove that a convex body (even highly symmetric) can have large volume ratio (close to the volume ratio of the cube of the same volume) but small minimal surface area (close to the surface area of the ball of the same volume). Notice that it is well known  that the opposite cannot happen (see (\ref{eq-minimal-surface-vr}) below).

Bourgain's slicing problem \cite{Bou1} asks whether any convex body in  $\mathbb{R}^n$ of volume $1$  has a hyperplane section of volume greater than $c>0$, where $c$ is an absolute constant. It follows from the work of the first named author \cite[Theorem 1]{Ko2} (Theorem \ref{thm-Ko} in Section 5) that if a centrally symmetric convex body $K$ has bounded outer volume ratio with respect to the class of intersection bodies (see Section 5 for definition), then the slicing problem has an affirmative answer for $K$ (actually this result extends to general measures in place of volume). In Section 5, we extend this result (however, not bodywise) as follows (Theorem \ref{thm-slicing}): If every convex body in $\mathbb{R}^n $ has a projection of dimension at least proportional to $n$ which is close to an intersection body (in the sense of outer volume ratio), then the slicing problem will have an affirmative answer. 
\smallbreak

\noindent {\bf Acknowledgment}. We are indebted to Matthieu Fradelizi for many valuable discussions  and anonymous referee for many suggestions which led to a better
presentation of the above results.

\section{Background and notation}
\hspace*{1.5em}We use the notation $a\sim b$ to declare that the ratio $a/b$ is bounded from above and from below by absolute constants. We denote by $D_n$ the standard $n$-dimensional Euclidean ball of volume $1,$ and by $B_2^n$ the Euclidean ball of radius 1. We denote the volume of $B_2^n$ by $\omega_n$.

A set $L$  is called a star body if it has non-empty interior, it is star-shaped at the origin and its radial function $\rho_L$ is continuous.  We remind that the radial function $\rho_L$ of $L$ is defined as:
$$\rho_L(u)=\max\{\lambda>0:\lambda u\in L\}\ ,\qquad u\in S^{n-1}\ ,$$
where $S^{n-1}=\{x\in\mathbb{R}^n:|x|=1\}$ is the unit sphere in $\mathbb{R}^n$.

In this section, $K$ will always denote a convex body in $\mathbb{R}^n$. The support function of $K$ is defined as
$$h_K(x)=\max_{y\in K}\langle x,y\rangle\ ,\qquad x\in\mathbb{R}^n \ .$$
The surface area measure (or curvature measure) $S_K$ of $K$ is a measure on $S^{n-1}$, defined by
$$S_K(\Omega)=\mathcal{H}^{n-1}\Big(\big\{x\in bdK:\exists u\in \Omega, \textnormal{ such that }\langle x,u\rangle=h_K(u)\big\}\Big)\ ,\qquad \Omega\textnormal{ is a Borel subset of }S^{n-1}\ ,$$
where $\mathcal{H}^{n-1}$ is the $(n-1)$-dimensional Haussdorff measure. If $S_K$ is absolutely continuous with respect to the Lebesgue measure, its density is denoted by $f_K$ and it is called the curvature function of $K$.

A symmetric convex body  $\Pi K$ is ``the projection body of $K$'' if its support function is given by 
\begin{equation}\label{eq:zonoid}
h_{\Pi K}(u)=|K|u^{\perp}|= \frac{1}{2}\int\limits_{S^{n-1}} |x \cdot u| d S_K(x) \mbox{  for all  } u\in S^{n-1}, 
\end{equation}
where the last equality follows from the Cauchy formula for the volume of the orthogonal projection. Next, if $K_1,\dots,K_n$ are compact convex sets in $\mathbb{R}^n$, we denote their mixed volume by $V(K_1,\dots,K_n)$. We refer to \cite{S} or \cite{Ga} for an extensive discussion on the theory of mixed volumes. Note that if $S(K)$ is the surface area of $K$, then
\begin{equation*}\label{eq-S(K)=V(...)}
S(K)=nV(K[n-1],B_2^n)=\int_{S^{n-1}}dS_K.
\end{equation*}
In general, $$V(K[n-1],L)=\frac{1}{n}\int_{S^{n-1}}h_LdS_K$$
and 
$$|K|=\frac{1}{n}\int_{S^{n-1}}h_KdS_K.$$ The Minkowski inequality for mixed volumes states that
\begin{equation}\label{eq-Minkowski}
V(K[n-1],L)\geq |K|^{(n-1)/n}|L|^{1/n},
\end{equation}
with equality if and only if $K$ and $L$ are homothetic.

As F. John \cite{J} (see also \cite{AGM}, page 50) proved, there exists a unique ellipsoid $JK$ of maximal volume contained in $K$, the so-called ``John ellipsoid'' of $K$. Since,  for $T\in GL(n)$ we have that $J(TK)=T(JK)$, there always exists $T\in SL(n)$, such that $J(TK)$ is a ball. The quantity $vr(K):=\big(|K|/|JK|\big)^{1/n}$ is called the volume ratio of $K$. 

Consider a convex body $K$, with $|K|=1$ and such that  $JK$ is an Euclidean  ball. Then,
$$1=|K|=\frac{1}{n}\int_{S^{n-1}}h_KdS_K\geq \frac{1}{n}\int_{S^{n-1}}h_{JK}dS_K=\frac{|JK|^{1/n}}{n\omega_n^{1/n}}S(K).$$
We know by Stirling's formula that $\omega_n^{1/n}\sim 1/\sqrt{n}$, hence
\begin{equation}\label{eq-surface-vr}
 S(K)\leq c\sqrt{n}\frac{1}{|JK|^{1/n}}=c\sqrt{n}vr(K)\ .
\end{equation}
Let us also define the quantity $$\partial(K):=\min_{T\in GL(n)}\frac{S(TK)}{|T K|^{(n-1)/n}}\ .$$
We say that $K$ is in minimal surface area position if $S(K)=\partial (K)$ and $|K|=1$ (see \cite{AGM}, Section 2.3). With this definition, (\ref{eq-surface-vr}) gives:
\begin{equation}\label{eq-minimal-surface-vr}
\partial(K)\leq c\sqrt{n}vr(K)\ . 
\end{equation}
 The parameters $\partial (K)$ and $vr(K)$ are affine invariants. A useful characterization of the minimal surface area position due to Petty \cite{Pe} states that a convex body $K$ of volume $1$ is in minimal surface area position if and only if its  surface area measure $S_K$ is isotropic, i.e.
\begin{equation}\label{eq-S_K-isotropic}
\int_{S^{n-1}}|\langle x,u\rangle|^2dS_K(x)=\frac{1}{n}S(K)\ ,\qquad \textnormal{ for all }u\in S^{n-1}\ .
\end{equation}

We say that the convex body $K$ of volume $1$ is in minimal mean width position (see \cite[Section 2.2]{AGM}) if
$$
\frac{1}{|S^{n-1}|}\int_{S^{n-1}} h_K(u) du \le  \frac{1}{|S^{n-1}|}\int_{S^{n-1}} h_{TK}(u) du 
$$
for every $T\in SL(n)$. A very useful result  which follows from  Figiel,  Tomczak-Jaegermann, Lewis and Pisier estimates (see  \cite[Corollary 6.5.3]{AGM})  gives that for a symmetric convex body $K\subset \R^n$ in minimal mean width position and of volume $1$ we get
\begin{equation}\label{eq:min}
\frac{1}{|S^{n-1}|}\int_{S^{n-1}} h_K(u) du \le C \sqrt{n} \log n,
\end{equation}
for some absolute constant $C>0$.

If $K$ contains the origin in its interior, the polar body $K^{\circ}$ of $K$ is defined to be the convex body 
$$K^{\circ}=\{x\in\mathbb{R}^n:\langle x,y\rangle\leq 1,\ \forall y\in K\}\ .$$
If $T\in GL(n)$ and $H\in G_{n,k}$, the following formulas hold:
$$(TK)^{\circ}=T^{-\ast}K^{\circ}\qquad \textnormal{and}\qquad (K|H)^{\circ}=K^{\circ}\cap H\ .$$
Here, $G_{n,k}$ denotes the Grassmannian manifold of $k$-dimensional subspaces of $\mathbb{R}^n$ and the notation $A|H$ denotes the orthogonal projection of  $A$ onto the subspace $H$.

Let us assume  that the origin is the centroid of $K$. The Blaschke-Santal\'o inequality (see \cite{San, RZ}) together with its reverse (see \cite{BM, RZ}) give:
\begin{equation}\label{eq-Santalo}
(|K|\cdot|K^{\circ}|)^{1/n}\sim 1/n\ .                                                                                                                                                                                                \end{equation}

Set $N(K,D_n)$ to be the covering number of $K$ by translates of $D_n$, i.e.
$$
N(K,D_n)=\min\{N\in\mathbb{N}:  \exists  \,\,  x_1,\dots,x_N\in\mathbb{R}^n,\textnormal{ such that }K\subseteq \cup_{i=1}^N(x_i+D_n)\}\ .
$$
Milman (\cite{M}, see also \cite[Chapter 8]{AGM}) proved that there exists an absolute constant $C>0$ and a linear image $K'$ of $K$ of volume 1, such that
$$\max\Big\{N(K',D_n), N(|K'^{\circ}|^{-1/n}K'^{\circ},D_n), |K'\cap D_n|^{-1},|(|K'^{\circ}|^{-1/n}K'^{\circ})\cap D_n|^{-1}\Big\}\leq C^{n}\ .$$
If the previous inequality holds for $K'$, we say that $K'$ is in $M$-position.

We will also need to use the notion of isotropic position of a convex body (see \cite{MP, BGVV}): there exists $T\in SL(n)$, such that
$$\int_{TK}\langle x,y\rangle^2dy=L_K^2|K|^{(n+2)/n}|x|^2\ ,\qquad x\in\mathbb{R}^n \,$$   we will call $TK$ an isotropic convex body or a body  in isotropic position. The parameter $L_K$ is called the isotropic constant of $K$. $L_K$ depends only on $K$ and it is invariant under invertible linear maps. It turns out that
\begin{equation}\label{eq-L_K-variational}
L_K^2=\frac{1}{n}\min_{T\in SL(n)}\frac{1}{|K|^{(n+2)/n}}\int_{TK}|x|^2dx=\frac{1}{n}\min_{T\in GL(n)}\frac{1}{|TK|^{(n+2)/n}}\int_{TK}|x|^2dx\ .
\end{equation}
It is a major problem to show that $L_K$ is uniformly bounded from above by an absolute constant (the fact that $L_K>c$ can be proved by comparison with an Euclidean ball; see  \cite{MP, BGVV}). The best estimate currently is due to Klartag \cite{K}: $L_K\leq C'n^{1/4}$ who removed the logarithmic term in the previous estimate of Bourgain \cite{Bou2}. It should be noted that (see \cite{MP, BGVV}): 
\begin{equation}\label{eq-hensley}
|K\cap u^{\perp}|\sim \bigg(\int_K\langle x,u\rangle^2dx\bigg)^{-1/2}\ ,\qquad u\in S^{n-1}\ . 
\end{equation}
Thus, the problem of bounding $L_K$ is equivalent to Bourgain's slicing problem.
\section{A remark on the reverse Loomis-Whitney inequality}
The main goal of this section is to provide a sharp asymptotic estimate for the quantity $\Lambda_n^{1/n}$ in the reverse Loomis-Whitney inequality.
\begin{theorem}\label{thm-reverse-lw}
There exists an absolute constant $c>0$, such that 
\begin{equation}\label{eq-thm-lw-1}
\Lambda_n\leq (c\sqrt{n})^n. 
\end{equation}
Moreover, there exists a symmetric convex body $L$ in $\mathbb{R}^n$, such that 
\begin{equation}\label{eq-thm-lw-2}\min_{\{e_1,\dots,e_n\}\in{\mathcal F}_n}\prod_{i=1}^n\big|L|e_i^{\perp}\big|\geq (c'\sqrt{n})^n|L|^{n-1}\ ,
\end{equation}
where $c'>0$ is another absolute constant. 
\end{theorem}
The proof of Theorem \ref{thm-reverse-lw} follows from two theorems due to K. Ball. The first is K. Ball's reverse isoperimetric inequality (for symmetric convex bodies):
\begin{oldtheorem}\label{oldthm-ball-1}\cite{Ba2}
Let $K$ be a convex body and $\Delta$ be a simplex in $\mathbb{R}^n$. Then,
$$\partial (K)\leq \partial (\Delta)\leq cn\ .$$
\end{oldtheorem}
The second is a remarkable example of a convex body whose projections all have large volumes (comparing to its volume):
\begin{oldtheorem}\label{oldthm-ball-2}\cite{Ba1}
There exists a symmetric convex body $L$ in $\mathbb{R}^n$ of volume 1, such that $\big|L|u^{\perp}\big|\sim\sqrt{n}$, for all $u\in S^{n-1}$. 
\end{oldtheorem}
\textit{}\\
\noindent{\it Proof of Theorem \ref{thm-reverse-lw}:} 
First, note that the convex body $L$ from Theorem \ref{oldthm-ball-2} serves as an example of a convex body that satisfies (\ref{eq-thm-lw-2}). 

To prove (\ref{eq-thm-lw-1}), we  may assume that $K$ is of volume 1. Let us first consider a convex body $M$ in minimal surface area position (thus $|M|=1$). It was observed  in \cite{GP}, that, in this case, the measure $S_M$ is isotropic. By equation (\ref{eq-S_K-isotropic}), Theorem \ref{oldthm-ball-1} and the Cauchy-Schwartz inequality 
\begin{eqnarray}\label{eq-proj-up-bound}
 \big|M|u^{\perp}\big|&=&\frac{1}{2}\int_{S^{n-1}}|\langle x,u\rangle|dS_M(x)\leq \frac{1}{2}\bigg(\int_{S^{n-1}}|\langle x,u\rangle|^2dS_M(x)\bigg)^{1/2} \bigg(\int_{S^{n-1}}dS_M(x)\bigg)^{1/2}\nonumber\\
 &=&\frac{1}{2}\bigg(\frac{1}{n}\int_{S^{n-1}}|x|^2dS_M(x)\bigg)^{1/2} \bigg(\int_{S^{n-1}}dS_M(x)\bigg)^{1/2}\nonumber\\
 &=&\frac{1}{2\sqrt{n}}\partial(M)|M|^{(n-1)/n}\leq c\sqrt{n}\ ,
\end{eqnarray}
for all $u\in S^{n-1}$.

Note that for every convex body $K$, $|K|=1$ there exists $T\in SL(n)$, such that $TK$ is in minimal surface area position. Since both (\ref{eq-S_K-isotropic}) and (\ref{eq-thm-lw-1}) are invariant under orthogonal transformations, we may  assume that $T$ is a diagonal matrix  with strictly positive entries. Then, $T^{-1}=\textnormal{diag}(\lambda_1,\dots,\lambda_n)$ for some $\lambda_1,\dots,\lambda_n>0$ with $\prod\lambda_i=1$. Consequently, if $\{e_1,\dots,e_n\}$ is the standard orthonormal basis in $\mathbb{R}^n$, it follows by (\ref{eq-proj-up-bound}) and the fact that $T^{-1}$ is diagonal that:
\begin{eqnarray*}
\prod_{i=1}^n\big|K|e_i^{\perp}\big|&=&\prod_{i=1}^n \big|T^{-1}T(K|e_i^{\perp})\big|=\prod_{i=1}^n\big|TK|e_i^{\perp}\big|\prod_{j\in\{1,\dots,n\}\setminus\{i\}}\lambda_j\\
&=&\prod_{i=1}^n\big|TK|e_i^{\perp}\big|\frac{1}{\lambda_i}\leq (c\sqrt{n})^n\prod_{i=1}^n\frac{1}{\lambda_i}=(c\sqrt{n})^n\ .
\end{eqnarray*}
\endproof

\section{Curvature measures and comparison of volumes}

Let us continue in another direction with the following observation:
\begin{proposition}\label{planar}
Let $K$, $L$ be centrally symmetric convex bodies in $\mathbb{R}^2$. If $S_K\leq S_L$, then $K\subseteq L$. 
\end{proposition}

\noindent{\it Proof:} Using the definition of projection body of a convex body (see equation (\ref{eq:zonoid})) we get that  $\Pi K\subseteq \Pi L$. However one can notice that for any symmetric convex body $M \subset \mathbb{R}^2$ we have $\Pi M=2O M$, where $O$ is the rotation by $\pi/2$. The result follows. \endproof
\begin{remark}
After a suitable translation, Proposition \ref{planar} remains true in the non-symmetric planar case as well. This is due to the additivity of the curvature measures of planar convex bodies. Indeed, set $\mu:=S_L-S_K$. Minkowski's Existence and Uniqueness Theorem (see \cite[Section 8.2]{S}) states that there exists a unique (up to translation) compact convex set $M$, such that $S_M=\mu$. Assuming that $K$, $L$ and $M$ contain the origin, we get: $L=K+M\supseteq K$, where we used the fact (see \cite[Section 8.3]{S}) that in the plane $S_{K+M}=S_K+S_M$ and again Minkowski's Existence and Uniqueness Theorem.  

\end{remark}
One cannot naturally expect Proposition \ref{planar} to hold true in any dimension. Indeed, let  $K=[-1,1]^{n-3}\times[-1/10,1/10]^2\times[-2,2]$ and $L=[-1,1]^n$, $n\ge 3$, then $S_K \le S_L$, but $K \not\subset L$.  From another point of view, one can always guarantee the comparison of volumes:
\begin{proposition}\label{notplanar}
Let $K$, $L$ be centrally symmetric convex bodies in $\mathbb{R}^n$. If $S_K\leq S_L$, then $|K| \le |L|$. 
\end{proposition}

\noindent{\it Proof:} Indeed, from $S_K\leq S_L$ we get $\int_{S^{n-1}} h_L S_K \le \int_{S^{n-1}} h_L S_L$, and can use  Minkowski inequality (\ref{eq-Minkowski}) to finish the proof.\endproof
Thus it is interesting to ask for comparison with a Euclidean Ball:
If $K$ is a symmetric convex body of volume $1$, such that $JK$ is an Euclidean ball and $S_K\geq S_{cD_n}$, for some absolute constant $c>0$, is it true that $K$ has bounded volume ratio? Actually, one might ask a weaker (as (\ref{eq-minimal-surface-vr}) shows) version of the previous question: 
\begin{question}\label{question 1}
Let $K$ be a symmetric convex body, such that $JK$ is a ball. If $S_K\geq S_{cD_n}$, is it true that the quantity $\partial(K)/\sqrt{n}$ is bounded from above by an absolute constant?
\end{question}
An extra motivation for Question \ref{question 1} is Proposition \ref{prop-conditional} (see below). As we show the answer to this question is negative. Recall that a convex body is called 1-symmetric if its symmetry group contains the symmetry group of the standard $n$-cube.  
\begin{theorem}\label{thm-curvature-vr}
There exists an absolute contant $c>0$, such that for each positive integer $n$, there exists an 1-symmetric convex body $L$ of volume 1  in $\mathbb{R}^n$, such that $S_L>cS_{D_n}$, but $\partial(L)/\sqrt{n}\geq c\sqrt{n}$.  
\end{theorem}
\begin{remark}\textit{}
\begin{enumerate}
\item Note that if $L$ as in Theorem \ref{thm-curvature-vr}, i.e. $1$-symmetric,  then  $JL$ is an Euclidean ball and $\partial(L)=S(L)$.
\item K. Ball's reverse isoperimetric inequality shows that the factor $1/\sqrt{n}$ gives the worst possible order in Theorem \ref{thm-curvature-vr}.
\item We also note that $cS_{D_n}=S_{c^{1/(n-1)}D_n}$, so the assumption in Theorem \ref{thm-curvature-vr} is much stronger than the assumption in Question \ref{question 1}.
\end{enumerate}
\end{remark}
\noindent{\it Proof of Theorem \ref{thm-curvature-vr}:}  Let $C_n$ be the standard $n$-cube of volume 1. Consider the convex body $K$ defined as the Blaschke average of $C_n$ and $D_n$:
$$S_K=\frac{1}{2}S_{C_n}+\frac{1}{2}S_{D_n}\ .$$
Note that $C_n$ is in $M$-position (this, follows, for example from covering $C_n$  by copies of $D_n$, see Lemma 8.1.3 in \cite{AGM}), thus:$$|C_n\cap D_n|\geq C^{-n} \ .$$
Using the Minkowski inequality (\ref{eq-Minkowski}), we obtain:
\begin{eqnarray*}
C^{-1}|K|^{(n-1)/n}&\leq& |C_n\cap D_n|^{1/n}|K|^{(n-1)/n}\leq\frac{1}{n}\int_{S^{n-1}}h_{C_n\cap D_n}dS_K\\
&=&\frac{1}{2n}\int_{S^{n-1}}h_{C_n\cap D_n}dS_{C_n}+\frac{1}{2n}\int_{S^{n-1}}h_{C_n\cap D_n}dS_{D_n}\\
&\leq&\frac{1}{2n}\int_{S^{n-1}}h_{C_n}dS_{C_n}+\frac{1}{2n}\int_{S^{n-1}}h_{ D_n}dS_{D_n}\\
&=&\frac{1}{2}|C_n|+\frac{1}{2}|D_n|=1\ .
\end{eqnarray*}
This shows that 
\begin{equation}\label{eq-|K|^{(n-1)/n}}
|K| ^{(n-1)/n}\leq C \ .
\end{equation}
Moreover, using the fact that $S(D_n)\sim \sqrt{n}$ and $S(C_n)=2n$, we get:
$$\frac{1}{2}\frac{2n+c'\sqrt{n}}{C}\leq \frac{S(C_n)+S(D_n)}{2C}\leq \frac{S(K)}{|K|^{(n-1)/n}}\ ,$$
therefore 
$$\frac{S(K)}{|K|^{(n-1)/n}}\geq c''n\ .$$
Set $L:=(1/|K|^{1/n})K$. Then, $|L|=1$ and  $$S_L\geq \frac{1}{2}\frac{1}{|K|^{(n-1)/n}}S_{D_n}\geq [1/(2C)]S_{D_n}\ .$$However, $$\partial(L)=\frac{S(K)}{|K|^{(n-1)/n}}\geq c''n\ ,$$as claimed. \endproof

Let $K_0$ be a centrally symmetric star body in $\mathbb{R}^n$. Let us recall the definition of the curvature image $C(K_0)$ body of $K_0$  which is defined via the curvature function of $C(K_0)$ (see \cite{L1} for more information):
$$f_{C(K_0)}(x)=\frac{\rho_{K_0}^{n+1}(x)}{n+1}\ ,\qquad x\in S^{n-1}\ .$$
A convex body $L$ is called ``body of elliptic type'' if there exists a convex body $K_0$, such that $L=C(K_0)$. 

\begin{question}\label{question 1'}
Is Question \ref{question 1} true for symmetric bodies of elliptic type? 
\end{question}
\begin{remark}\label{prop-conditional} We note that  following   ideas of  the proof from \cite[Proposition 5.3]{Sa} one can show that 
if Question \ref{question 1'} had an affirmative answer, then the slicing problem would have an affirmative answer as well.
\end{remark}
As shown in \cite{Sa}, for any symmetric convex body $K$, $L_K\sim \partial(C(K))$ (so the reverse implication to the statement from Remark \ref{prop-conditional} obviously holds). Hence, it follows by (\ref{eq-minimal-surface-vr}) that if all centrally symmetric bodies of elliptic type have uniformly bounded volume ratios, then the isotropic constant would be uniformly bounded. It would be, therefore, natural to ask if the reverse inequality of (\ref{eq-minimal-surface-vr}) is true.  

\begin{question}\label{question 2}
Is there an absolute constant $c'>0$. such that 
$$c'\sqrt{n}vr(K)\leq \partial(K)\ ,$$
for all centrally symmetric bodies of elliptic type $K$?
\end{question}
If Question \ref{question 2} had an affirmative answer, then for any symmetric convex body $K$, the following would hold true: $L_K$ is bounded if and only if $vr(C(K))$ is bounded. 
This provides a good motivation for our next result: a reverse inequality to (\ref{eq-minimal-surface-vr}) cannot hold true for general symmetric convex bodies.
\begin{theorem}\label{counterexample}
There exists a convex body $K$ in $\mathbb{R}^n$, such that $$\partial(K)\leq c'vr(K)\sqrt{\log(n+1)}\ .$$  
\end{theorem}
\begin{remark}\text{}
\begin{enumerate}
\item Actually $K$, from Theorem \ref{counterexample} can be taken to be 1-symmetric.
\item It follows by K. Ball's volume ratio inequality \cite{Ba2} and the classical isoperimetric inequality that Theorem \ref{counterexample} gives the worst possible case up to a logarithmic factor. 
\end{enumerate}
\end{remark}
\noindent{\it Proof of Theorem \ref{counterexample}:} We first notice that the mean width of  $B_{\infty}^n$ (the unit ball of $\ell_\infty^n$) is of the order $\sqrt{(\log (n+1))/n}$ (this can be calculated by passing to the integral over the gaussian measure and using standard estimates for the sequence of standard normal variables, see \cite{AGM} and \cite[Chapter 3, page 79]{LT} for more details). Therefore, integrating in polar coordinates, we get:
\begin{eqnarray*}
\int_{D_n}\|x\|_{\infty}dx&=&\frac{1}{n+1}\int_{S^{n-1}}\|x\|_{\infty}\rho^{n+1}_{D_n}dx\\
&=&\frac{\omega_n^{-\frac{n+1}{n}}}{n+1}|S^{n-1}|\int_{S^{n-1}}\|x\|_{\infty}\frac{dx}{|S^{n-1}|}\\
&\sim&\omega_n^{-1/n}\int_{S^{n-1}}\|x\|_{\infty}d\sigma(x)\\
&\sim&\sqrt{n}\sqrt{(\log (n+1))/n}=\sqrt{\log(n+1)}\ ,
\end{eqnarray*}
where $\sigma$ is the Haar probability measure on $S^{n-1}$. We note that 
\begin{eqnarray*}
\int_{D_n}\|x\|_{\infty}dx=\int_0^{\infty}|D_n\cap \{\|x\|_{\infty}> s\}|ds=\int_0^{\infty}|D_n\setminus sB^n_{\infty}|ds \ .
\end{eqnarray*}
Let $s_0>0$ be such that 
\begin{equation}\label{eq-s_0}
|D_n\cap s_0B_{\infty}^n|=|D_n\setminus s_0B^n_{\infty}|=1/2\ . 
\end{equation}

\begin{eqnarray*}
 \int_{D_n}\|x\|_{\infty}dx&=&\int_0^{s_0}|D_n\setminus sB^n_{\infty}|ds+\int_{s_0}^{\infty}|D_n\setminus sB^n_{\infty}|ds\\
 &>&\int_0^{s_0}|D_n\setminus s_0B^n_{\infty}|ds+\int_{s_0}^{\infty}|D_n\setminus sB^n_{\infty}|ds>s_0/2\ .
\end{eqnarray*}
It follows that $s_0\leq C_1\sqrt{\log(n+1)}$. Take, now, $$K:=D_n\cap s_0B_{\infty}^n\ .$$
Note that $K$ is 1-symmetric, with  $|K|=1/2$ and  thus
$$\partial(K)=S(K)\leq S(D_n)\leq C_2\sqrt{n}\ .$$
Moreover, since $K$ is 1-symmetric, $JK$ is an Euclidean ball and the largest Euclidean ball contained in $K$. Using the definition of $K$ we get that $JK$  is not larger then the largest Euclidean  ball contained in $s_0B_{\infty}^n$, so $$JK\subseteq s_0B_2^n\leq (1/c_1)\sqrt{\log(n+1)}B_2^n\ ,$$for some absolute constant $c_1>0$. Thus, $$vr(K)\geq c_1\bigg(\frac{|K|}{|\sqrt{\log(n+1)}B_2^n|}\bigg)^{1/n}\geq c_1\frac{1/2}{\sqrt{\log(n+1)}\omega_n^{1/n}}\geq c_2\frac{\sqrt{n}}{\sqrt{\log(n+1)}}\geq \frac{c_2}{C_2}\frac{\partial(K)}{\sqrt{\log(n+1)}}\ .$$ \endproof
\begin{remark}
By a well known lemma of C. Borell \cite{Bor}, we have
$$|D_n\setminus sB^n_{\infty}|\leq C_0e^{-c_0s}\ ,$$
for $s\geq 2s_0$, where $c_0, C_0>0$ are absolute constants and $s_0$ is defined by (\ref{eq-s_0}). Therefore,
$$\int_{2s_0}^{\infty}|D_n\setminus sB_{\infty}^n|ds\leq \int_0^{\infty}C_0e^{-c_0s}ds\leq C_0'$$
and
$$\int_0^{2s_0}|D_n\setminus sB^n_{\infty}ds|\leq 2s_0|D_n|=2s_0\ .$$
Thus, $s_0\geq c_3\sqrt{\log(n+1)}$, which shows that the logarithmic factor in the example of Theorem \ref{counterexample} cannot be removed. 
\end{remark}
\begin{remark}
 Note that if $K$ is the example from Theorem \ref{counterexample} (recall that $|K|=1/2$), by (\ref{eq-S_K-isotropic}) and (\ref{eq-proj-up-bound}), we have
$ |K|u^{\perp}|\leq c_4\ ,$ where we used the fact that since $K$ is 1-symmetric, $S_K$ is isotropic. Also, $|K|u^{\perp}|\geq |K\cap u^{\perp}|\geq c_{5}|K|^{(n-1)/n}\geq c_6$, since $K$ is in isotropic position and it has bounded isotropic constant, as 1-symmetric. Consequently (after a suitable dilation), we may take $K$ in Theorem \ref{counterexample} such that all its projections have volumes of constant order, when it is in minimal surface area position. 
\end{remark}

The next two theorems continue the discussion started in Proposition \ref{notplanar}.  Our goal is to  provide bounds on the difference of volumes of convex bodies via its curvature functions. 
\begin{thm}\label{th:min} Let $K, L \subset \R^n$ be  convex bodies with absolutely continuous surface area measure, such that $f_K(\theta) \le f_L(\theta)$, for all $\theta\in S^{n-1}$. Then 
$$
|L|^{\frac{n-1}{n}} - |K|^{\frac{n-1}{n}}  \ge \omega_n^{\frac{n-1}{n}} \min\limits_{\theta \in S^{n-1}} (f_L(\theta)- f_K(\theta)).
$$
\end{thm}
\noindent{\it Proof:}  Consider two convex bodies $K,L \subset \R^n$ such that 
\begin{equation}\label{eq:minf}
f_K(\theta) \le f_L(\theta)- \epsilon,
\end{equation}
 for some $\epsilon \ge 0$ and all $\theta \in S^{n-1}$. Then
$$
\int\limits_{S^{n-1}} h_L(\theta) f_K(\theta) d \theta \le \int\limits_{S^{n-1}} h_L(\theta) f_L(\theta) d \theta - \epsilon \int\limits_{S^{n-1}} h_L(\theta) d \theta
$$
or
$$
V(K[n-1], L) \le |L| -\epsilon V(B_2^n[n-1], L).
$$
Applying the Minkowski inequality (\ref{eq-Minkowski}) to $V(K[n-1], L)$ and $V(B_2^n[n-1], L)$ we get
$$
|L|^{\frac{1}{n}}|K|^{\frac{n-1}{n}} \le |L| -\epsilon |L|^{\frac{1}{n}}|B_2^n|^{\frac{n-1}{n}}.
$$
Thus
$$
|L|^{\frac{n-1}{n}} - |K|^{\frac{n-1}{n}} \ge \epsilon |B_2^n|^{\frac{n-1}{n}}.
$$
To finish the proof of the theorem we note that (\ref{eq:minf}) is always true with $\epsilon = \min\limits_{\theta \in S^{n-1}} (f_L(\theta)- f_K(\theta))$.
\endproof
We note that if we can consider $K$ to be $rB_2^n$ in the above theorem and send $r \to 0$ to get that 
$$
|L|^{\frac{n-1}{n}}  \ge \omega_n^{\frac{n-1}{n}} \min\limits_{\theta \in S^{n-1}} (f_L(\theta)),
$$
for any $L$ with absolutely continuous surface area measure.

\begin{thm}\label{th:max} There exists an absolute constant $C>0$ such that for any $K, L \subset \R^n$, two convex bodies with absolutely continuous surface area measure, with $f_L(\theta) \ge f_K(\theta),$ for all $\theta \in S^{n-1}$ and such that $K$ is in the minimal width position,  we have
$$
|L|^{\frac{n-1}{n}} - |K|^{\frac{n-1}{n}}  \le C \log n  \,\, \omega_n^{\frac{n-1}{n}} \max \limits_{\theta \in S^{n-1}} (f_L(\theta)- f_K(\theta)).
$$
\end{thm}
\noindent{\it Proof:}  Consider two convex bodies $K,L \subset \R^n$ such that 
\begin{equation}\label{eq:max}
f_L(\theta) \le f_K(\theta)+ \epsilon,
\end{equation}
 for some $\epsilon \ge 0$ and all $\theta \in S^{n-1}$. Then
$$
\int\limits_{S^{n-1}} h_K(\theta) f_L(\theta) d \theta \le \int\limits_{S^{n-1}} h_K(\theta) f_K(\theta) d \theta + \epsilon \int\limits_{S^{n-1}} h_K(\theta) d \theta
$$
or
$$
V(L[n-1], K) \le |K| +\epsilon \frac{1}{n} \int\limits_{S^{n-1}} h_K(\theta) d \theta.
$$
Now we can apply  the Minkowski inequality (\ref{eq-Minkowski}) to $V(L[n-1], K)$.  We can also use the fact that $K$ is in minimal width position  and apply inequality (\ref{eq:min})  to claim that
$$
 \int\limits_{S^{n-1}} h_K(\theta) d \theta \le C |S^{n-1}| |K|^{\frac{1}{n}} \sqrt{n} \log n =  C \omega_n |K|^{\frac{1}{n}} n \sqrt{n} \log n.
$$
Thus
$$
|K|^{\frac{1}{n}}|L|^{\frac{n-1}{n}} \le |K| + \epsilon \frac{1}{n}  C \omega_n |K|^{\frac{1}{n}} n \sqrt{n} \log n.
$$
or
$$
|L|^{\frac{n-1}{n}} - |K|^{\frac{n-1}{n}} \le \epsilon C |B_2^n|^{\frac{n-1}{n}} \log n.
$$
To finish the proof of the theorem we note that (\ref{eq:max}) is always true for $\epsilon = \max\limits_{\theta \in S^{n-1}} (f_L(\theta)- f_K(\theta))$.
\endproof

\section{A note on the slicing problem}

Recall that the notion of an {\it intersection body of a star body} was
 introduced by E. Lutwak \cite{L2}:
$IL$ is called the intersection body of $L$ if the radial function of $IL$ in
every direction is equal to the $(n-1)$-dimensional volume of the
central hyperplane section of $L$ perpendicular to this direction:
$\forall \xi\in S^{n-1},$
$$
\label{eq:int}\rho_{IL}(u)= |L\cap u^{\perp}|,
$$
As shown in \cite{L2} (using a result of Busemann's theorem \cite{Bu}), if $L$ happens to be a symmetric convex body, then $IL$ is also a symmetric convex body.

A more general class of {\it intersection bodies} was defined by R.
Gardner \cite{Ga1} and G. Zhang \cite{Zh}
 as the closure of intersection bodies of star bodies in
the radial metric $d(K,L)=\sup_{\xi\in S^{n-1}} |\rho_K(\xi)-
\rho_L(\xi)|$.  We refer reader to books \cite{Ga, Ko1, KoY, RZ} for
more information on definition and properties of Intersection body
and their applications in Convex Geometry and Geometric Tomography.

Define the quantity $L_n:=\max\{L_K:K\textnormal{ is a convex body in }\mathbb{R}^n\}.$
For $n,k\in\mathbb{N}$, $k\leq n$ and $c>0$, define the class $\mathcal{C}_{n,k,c}$ of convex
bodies $K$ in $\mathbb{R}^n$, that have a $l$-dimensional projection $P$, with $l\geq k$, for which there exists an $l$-dimensional intersection body, such that $P\subseteq L$ and $|L|/|P|\leq c^n$. Finally, define $I_{n,k}$ to be the smallest constant $t$, for which $\mathcal{C}_{n,k,t}$ contains all centrally symmetric convex bodies in $\mathbb{R}^n$.

The following was proved in  \cite[Corollary 1]{Ko2}:
\begin{oldtheorem}\label{thm-Ko}
There exists an absolute constant $\overline{C}_0>0$, such that if $L$ is a symmetric convex body in $\mathbb{R}^n$, such that $L\in \mathcal{C}_{n,n,a}$, for some $a>0$ and $\mu $ is any even measure with continuous density in $\mathbb{R}^n$, then  
\begin{equation}\label{eq-general-slicing}
\mu(L)\leq (\overline{C}_0/a)\max_{u\in S^{n-1}}\mu(L\cap u^{\perp})|L|^{1/n} \ . 
\end{equation}
\end{oldtheorem}
It follows by Theorem \ref{thm-Ko},  applied to  $\mu$ being standard Lebesgue measure,  together with  (\ref{eq-hensley}) that $L_n\leq \overline{C}_0'I_{n,n}$, for some absolute constant $\overline{C}_0'>0$. Our goal is to replace $I_{n,n}$ with $I_{n,k}$, where $k$ is at least proportional to $n$.
\begin{theorem}\label{thm-slicing}
There exists an absolute constant $\widetilde{C}>0$, such that for any $n\in\mathbb{N}$ and for any $\lambda\in(0,1]$, it holds $$L_n\leq \widetilde{C}^{1/\lambda}I_{n,\lfloor\lambda n\rfloor}\ ,$$where $\lfloor\cdot\rfloor$ denotes the floor function.
\end{theorem}
We also note that if $K$ has bounded outer volume ratio (in which case  $K$ also has bounded isotropic constant), then  $K^\circ$ has a bounded volume ratio and thus (see for example \cite[Theorem 6.1]{Pi}) $K^\circ$ has an almost Euclidean  section of proportional dimension, which gives a projection of $K$  of proportional dimension which is almost Euclidean. The following question is therefore natural:
\begin{question}
Is it true that any convex body $K \subset \mathbb{R}^n$ has projections of at least proportional dimension that are close to some intersection body in the sense of the Banach-Mazur distance?
\end{question}

Before proving Theorem \ref{thm-slicing}, we will need some geometric statements.

\begin{lemma}\label{lemma-small-projection}
Let $K$ be a symmetric convex body of volume 1, in $M$ position. There exists an absolute constant $C_0>0$, such that for any subspace $H$ of $\mathbb{R}^n$, it holds
$$1/C_0\leq|K\cap H|^{1/n}\leq |K|H|^{1/n}\leq C_0\ .$$
\end{lemma}
\noindent{\it Proof:} By assumption, there exist points $x_1,\dots,x_N\in\mathbb{R}^n$, with $N<C^n$, such that $$K\subseteq \bigcup_{i=1}^N(x_i+D_n)\ .$$
Let $H\in G_{n,k}$. Then, since $D_n|H=\omega_n^{-1/n}B_2^k$, we obtain
\begin{eqnarray*}
|K|H|&\leq& \Big|\bigcup_{i=1}^N(x_i+D_n) | H \Big|\leq N|D_n|H|
= N\omega_n^{-k/n}|B_2^k|\\
&=&N\omega_n^{-k/n}\omega_k\leq NC'^k(\sqrt{n})^{k}(\sqrt{k})^{-k}
\leq   (CC')^n(\sqrt{n/k})^k \ ,
\end{eqnarray*}
where we used the fact that $\omega_n^{1/n}\sim 1/\sqrt{n}$. 
Thus, $$|K|H|^{1/n}\leq \overline{C} \ .$$
On the other hand, if we replace $K$ by $(1/|K^{\circ}|)^{1/n}K^{\circ}$ and since $\big|(1/|K^{\circ}|)^{1/n}K^{\circ}\big|=1$, we get:
$$\Big|\big(\frac{1}{|K^{\circ}|^{1/n}}K^{\circ}\big)|H\Big|^{1/n}\leq \overline{C} \ ,$$
which gives
\begin{equation}\label{eq-n^2}
|K^{\circ}|H|^{1/n}\leq \overline{C}|K^{\circ}|^{k/n^2}\ .
\end{equation}
Using (\ref{eq-Santalo}) and the fact that $(K^{\circ}|H)^{\circ}=K\cap H$, we obtain
$$(c_1/k)^{k/n}|K\cap H|^{-1/n}\leq |K^{\circ}|H|^{1/n}.$$
This together with (\ref{eq-n^2}) give
$$(c_1/k)^{k/n}|K\cap H|^{-1/n}\leq \overline{C}|K^{\circ}|^{k/n^2}\leq \overline{C}(c_2/n)^{k/n}$$
or $$|K\cap H|^{1/n}\geq (nc_1/(c_2k))^{k/n}/\overline{C}\geq \overline{c}\ ,$$as required. \endproof
The next proposition is based on the existence of an $M$-position.
\begin{proposition}\label{prop-small-proj}
Let $K'$ be a symmetric convex body of volume 1, be such that $|K'\cap H|^{1/n}\leq C'$, for all subspaces $H$ of $\mathbb{R}^n$, for some absolute constant $C'>0$. Then, there exists an absolute constant $C'''$, such that $$|K'|H|^{1/n}\leq C'''\ ,$$for all subspaces $H$ of $\mathbb{R}^n$. 
\end{proposition}
\noindent{\it Proof:}  Consider a convex body  $K$ of volume $1$, in $M$-position, such that $K'=TK$, for some $T\in SL(n)$. As in the proof of Theorem \ref{thm-reverse-lw}, we may assume that $T$ is diagonal; actually we are only going to need the fact that $T$ is symmetric, i.e. $T^{\ast}=T$. We will make use of a special case of a result from \cite{F} (see also \cite[(5.28)]{S}), stating that if $M$ is a convex body in $\mathbb{R}^n$, $H\in G_{n,k}$ and $P$ is a convex body in $H^{\perp}$, then
\begin{equation}\label{eq-projections-general}
|M|H|=\binom{n}{k}\frac{1}{|P|}V\big(M[k],P[n-k]\big) \ . 
\end{equation}
Let $k\in\{1,\dots,n-1\}$ and $H\in G_{n,k}$.  Now we apply Lemma \ref{lemma-small-projection} together with (\ref{eq-Santalo}) to get:
\begin{equation}\label{eq-polar-sections}
\big|\big((TK)\cap H\big)^{\circ}\big|^{1/n}\geq (c_1/k)^{k/n}|(TK)\cap H|^{-1/n}\geq (c_1/k)^{k/n}/c' 
\end{equation}
and similarly,
\begin{equation}\label{eq-polar-projections}
 |K^{\circ}|H|^{1/n}\leq (c_2/k)^{k/n}|K\cap H|^{-1/n}\leq (c_2/k)^{k/n}/c''\ .
\end{equation}
Moreover, for any convex body $P$ in $H^{\perp}$, using (\ref{eq-projections-general}) we obtain:
\begin{eqnarray}\label{eq-after-polar-projections}
\big|(TK)^{\circ}|H\big|&=&\binom{n}{k}\frac{1}{|P|}V\big((T^{-\ast}K^{\circ})[k],P[n-k]\big)\nonumber \\
&=&\binom{n}{k}\frac{1}{|P|}V\big(K^{\circ}[k],(T^{\ast}P)[n-k]\big)\nonumber\\
&=&\binom{n}{k}\frac{1}{|P|}V\big(K^{\circ}[k],(TP)[n-k]\big)\nonumber\\
&=&\frac{|TP|}{|P|}\binom{n}{k}\frac{1}{|TP|}V\big(K^{\circ}[k],(TP)[n-k]\big)\nonumber\\
&=&\frac{|TP|}{|P|}\big|K^{\circ}|(T^{-\ast}H)^{\perp}\big|\ ,
\end{eqnarray}
where we used the fact that $TP$ is a convex body in $(T^{-\ast}H)^{\perp}$.
Thus, (\ref{eq-polar-sections}), (\ref{eq-polar-projections}) and (\ref{eq-after-polar-projections}) imply
$$\bigg(\frac{|TP|}{|P|}\bigg)^{1/n}\geq c''' \ .$$
Replacing $P$ by $T^{-1}P$, we get:
\begin{equation}\label{eq-T^{-1}P}
 |P|^{1/n}\geq c'''|T^{-1}P|^{1/n}\ ,
\end{equation}
for all convex bodies $P$ in $H$, for all subspaces $H$ of $\mathbb{R}^n$. On the other hand, as before one can write:
\begin{eqnarray*}
|(TK)|H|&=&\binom{n}{k}\frac{1}{|P|}V\big((TK)[k],P[n-k]\big)\\
&=&\binom{n}{k}\frac{1}{|P|}V\big(K[k],(T^{-1}P)[n-k]\big)\\
&=&\frac{|T^{-1}P|}{|P|}\binom{n}{k}\frac{1}{|T^{-1}P|}V\big(K[k],(T^{-1}P)[n-k]\big)\\
&=&\frac{|T^{-1}P|}{|P|}\big|K|(T^{\ast}H)\big|\ .
\end{eqnarray*}
Therefore, the previous identity, Lemma \ref{lemma-small-projection} and (\ref{eq-T^{-1}P}) give:
$$|K'|H|^{1/n}=\big|(TK)|H\big|^{1/n}\leq \big(|T^{-1}P|/|P|\big)^{1/n}\big|K|(T^{\ast}H)\big|^{1/n}\leq C'''/c'''\ ,$$
for any subspace $H$ of $\mathbb{R}^n$. \endproof
\begin{proposition}\label{prop-slicing-main}
Let $\lambda\in(0,1]$ and $K$ be a symmetric convex body in $\mathbb{R}^n$ of volume $1$, such that $\big|K|H\big|^{1/n}\leq A$, for some $A>0$ and for all subspaces $H$ of $\mathbb{R}^n$. If, in addition, $K\in \mathcal{C}_{n,k,t}$, for some $t>0$, $k\geq \lambda n$, then there exists an $(n-1)$-dimensional subspace $F$, such that $$|K\cap F|\geq \frac{\widetilde{c}}{A^{1/\lambda}t}\ ,$$ for some absolute constant $\widetilde{c}>0$.
\end{proposition}
 Let  $H\in G_{n,k}$. By Fubini's Theorem, we have:
\begin{equation}\label{eq-1st-fubini}
1=|K|=\int_{K|H}|K\cap (H^{\perp}+x)|dx\ .
\end{equation}
Let $G$ be a codimension $1$ subspace of $H$. Then, the subspace $F:=\textnormal{span}\{G\cup H^{\perp}\}$
has dimension $n-1$. The following claim is the key step in the proof of Proposition \ref{prop-slicing-main}:\\

\noindent\textbf{Claim.} $(K\cap F)|H=(K|H)\cap G$ and $(K\cap F)\cap (H^{\perp}+x)=K\cap (H^{\perp}+x)$, for all $x\in G$. 

\noindent{\it Proof of the claim:} To prove the first part, let $x\in(K|H)\cap G $. Since $x\in K|H$, there exists $y\in H^{\perp}$, such that $x+y\in K$. Since $x\in G$, $x+y\in F$, it follows that $x+y\in K\cap F$, so $x\in (K\cap F)|H$. Thus, $(K\cap F)|H\subseteq (K|H)\cap G$. Conversely, let $x\in (K\cap F)|H$. Then, there exists $y\in H^{\perp}$, such that $x+y\in K\cap F$. Thus, $x+y\in K$ and $x+y\in F$, so $x\in K|H$ and $x\in F|H=G$. Consequently, $x\in (K|H)\cap G$, which shows that $ (K|H)\cap G \subseteq(K\cap F)|H$, as required. To prove the second part, note that for any $x\in G$, $x+H^{\perp}\subseteq\textnormal{span}(G\cup H^{\perp})=F$, thus $K\cap(x+H^{\perp})\subseteq (K\cap F)\cap (x+H^{\perp})$. Since the opposite inclusion is trivial, our claim is proved. \endproof

\noindent{\it Proof of Proposition \ref{prop-slicing-main}:} Using again Fubini's Theorem and the previous claim, we get:
\begin{equation}\label{eq-2nd-fubini}
|K\cap F|=\int_{(K\cap F)|H}|K\cap F\cap (H^{\perp}+x)|dx=\int_{(K|H)\cap G}|K\cap F\cap (H^{\perp}+x)|dx=\int_{(K|H)\cap G}|K\cap (H^{\perp}+x)|dx\ .
\end{equation}
Assume, now, that $K|H\subseteq L$ and $|L|/|K|H|<t^n$, where $L$ is an intersection body in some $k$-dimensional subspace $H$ of $\mathbb{R}^n$. 

Set
$\mu$ to be the measure with density $H\in x\mapsto |K\cap (H^{\perp}+x)|$ in $H$ (i.e. $\mu$ is the marginal of the uniform measure on $K$ with respect to the subspace $H$; see e.g. \cite{Ba3}) and $L:=K|H\subseteq H$ . Using Theorem \ref{thm-Ko}, we conclude that there exists a codimension-1 subspace $G$ of $H$, such that:
\begin{equation}\label{eq-slicing-last}
\mu((K|H)\cap G)|K|H|^{1/k}\geq (\widetilde{c}/t)\mu(K|H)\ ,
\end{equation}
for some absolute constant $\widetilde{c}>0$. Since $$|K|H|^{1/k}\leq A^{1/\lambda}\ ,$$ if we set $F=\textnormal{span}\{G\cup H^{\perp}\}$, applying (\ref{eq-1st-fubini}) and (\ref{eq-2nd-fubini}) to (\ref{eq-slicing-last}),
we obtain:
$$|K\cap F|=\mu((K|H)\cap G)\geq \frac{\widetilde{c}_0}{A^{1/\lambda}t}\ ,$$as claimed. \endproof

\noindent{\it Proof of Theorem \ref{thm-slicing}.}  As shown in \cite{K}, there exists a symmetric convex body $K$, such that $L_K\geq c_0L_n$. We may assume that $K$ is of volume 1 and in isotropic position. Then, as it was proved in \cite[Corollary 3.5]{BKM}, $|K\cap H|^{1/n}$ is bounded from above by an absolute constant. It follows immediately by Propositions \ref{prop-small-proj}, \ref{prop-slicing-main} and the definition of $I_{n,k}$ that there exists an $(n-1)$-dimensional subspace $F$, such that 
$$|K\cap F|\geq \frac{\widetilde{c}}{C^{1/\lambda}I_{n,[\lambda n]}}\ ,$$
for some absolute constant $C>0$. Our claim follows immediately from the previous inequality and (\ref{eq-hensley}). \endproof

We remark that our method can be used to increase the class of symmetric convex bodies that are known to satisfy (\ref{eq-general-slicing}) (with an absolute constant). Indeed, one can replace the Lebesgue measure in Proposition \ref{prop-slicing-main} by any even measure (one has to repeat the steps of the proof) to obtain the following: 
\begin{proposition}\label{prop-main-general}
Let $\mu$ be an even measure in $\mathbb{R}^n$, with continuous density, $\lambda\in(0,1]$ and $K$ be a 0-symmetric convex body in $\mathbb{R}^n$, such that $\big|K|H\big|^{1/k}\leq A|K|^{1/n}$, for some $A>0$ and for all $k\geq \lambda n$ and $H\in G_{n,k}$. If, in addition, $K\in \mathcal{C}_{n,k_0,t}$, for some $t>0$ and for some $k_0\geq \lambda n$, then there exists an $(n-1)$-dimensional subspace $F$, such that $$\mu(K\cap F)|K|^{1/n}\geq \frac{\widetilde{c}}{A^{1/\lambda}t}\mu(K)\ ,$$ for some absolute constant $\widetilde{c}>0$. 
\end{proposition}
In particular, it follows from the previous statement and from Lemma \ref{lemma-small-projection} that if $K\in \mathcal{C}_{n,k,c}$, for some $k\geq \lambda n$ and $K$ is in $M$-position, then (\ref{eq-general-slicing}) holds for $K$ (with an absolute constant) for all even measures $\mu$ with continuous density. As an example, take $K=aB_2^n\times L$, where $L$ is any convex body in $M$-position, $\textnormal{dim}L=n$ and $|L|=|aB_2^n|$.

\end{document}